\documentclass[10pt,twoside]{article}
\usepackage{graphicx}
\usepackage{amsmath}
\usepackage{Latex-document}
\newtheorem{definition}{Definition}
\newtheorem{theorem}{Theorem}
\newtheorem{remark}{Remark}

\title{\bf Aging and Spin-glass Dynamics\vskip 6mm}

\author{G\'erard Ben Arous\vspace*{-0.5cm}\thanks{Department of Mathematics,
Ecole Polytechnique Federale de Lausanne, 1005 Lausanne, Switzerland. E-mail: gerard.benarous@epfl.ch}}

\date{\vspace{-8mm}}

\begin{document}

\maketitle

\thispagestyle{first} \setcounter{page}{3}

\begin{abstract}

\vskip 3mm

We survey the recent mathematical results about aging in certain simple disordered models. We start by the
Bouchaud trap model. We then survey the results obtained for simple models of spin-glass dynamics, like the REM
(the Random Energy Model, which is well approximated by the Bouchaud model on the complete graph), then the
spherical Sherrington-Kirkpatrick model. We will insist on the differences in phenomenology for different types of
aging in different time scales and different models. This talk is based on joint works with A.Bovier, J.Cerny,
A.Dembo, V.Gayrard, A.Guionnet, as well as works by C.Newman, R.Fontes, M.Isopi, D.Stein.

\vskip 4.5mm

\noindent {\bf 2000 Mathematics Subject Classification:} 60H10,
60K35, 82C44, 82C31, 20J05.

\noindent {\bf Keywords and Phrases:} Aging, Spin Glass, Random
Media, Trapping models, Statistical mechanics.
\end{abstract}

\vskip 12mm

\section{Introduction} \label{section 1}\setzero

\vskip-5mm \hspace{5mm}

Aging is an interesting long-time property of dynamics in complex disordered media, and in particular in certain
random media. A system ages when its decorrelation properties are age-dependent: the older the system the longer
it takes to forget its past. Aging has been heavily studied both experimentally, numerically and theoretically by
physicists, in particular in the context of spin-glass dynamics, but the mathematical litterature is still rather
sparse. Our interest in aging stemmed from the study of dynamics of mean field spin glasses models, and more
precisely of Langevin dynamics for the Sherrington Kirkpatrick model. This remains the ultimate goal, far from
being achievable. But we will survey some of the partial progress which has been made both in short time scales
(for the spherical SK model, see section 4) or much longer time scales (for the REM dynamics, see section 3). We
begin with aging the Bouchaud trap model on various graphs (see section 2)even though it is not directly related
to Spin Glasses, because the mechanism of aging there is close to the one for the REM. This talk is based on joint
works with A.Bovier, J.Cerny, A.Dembo, V.Gayrard, A.Guionnet, as well as works by C.Newman, R.Fontes, M.Isopi.

\section{Bouchaud's random trap model} \label{section 2}

\setzero\vskip-5mm \hspace{5mm}

This model is simple model of a random walk trapped (or rather
slowed down) by random wells. It is nevertheless quite rich. It
has been introduced by Bouchaud and coworkers (see [9] for a
beautiful general survey) as an ansatz to understand ``activated"
dynamics of spin glasses.

\subsection{The random energy landscape} \label{subsection 2.1}

\vskip-5mm \hspace{5mm}

Consider $G=(V,B)$ a graph where $V$ is the set of vertices and
$B$ the set of bonds. We introduce a ``random energy landscape" on
$V$, a collection $(E_x)_{x\in V}$ of i.i.d non-negative random
variables, indexed by the vertices of the graph and exponentially
distributed (with mean one) , i.e
$$ P(E_x > a)= e^{-a}. $$
From these random variables, we define a random measure $\tau$ on
$V$ by:
$$ \tau (x) = e^{\beta E_x}. $$
Here $\beta > 0 $ is an inverse temperature parameter. We will denote by $\alpha=\frac{1}{\beta}$
and concentrate on the low temperature phase, i.e $\alpha < 1$.

The physicists would see the set of vertices $V$ as a set of ``
favourable valleys" for the configurations of a much more complex
system (for instance a spin-glass), $ - E_x$ as the energy of the
bottom of the valley $x$. Since these energies are usually the
extreme values of some other random landscape (see below the
discussion about the REM), the hypothesis that they are
exponential random variables is reasonable. The sites $x$ where
$E_x$ is large are seen as `` very favourable valleys", or
equivalently as very deep traps, where the system should stay for
a long time, in any sensible definition of the dynamics.

\subsection{Bouchaud's random walk} \label{subsection 2.2}

\vskip-5mm \hspace{5mm}

We consider the continuous-time Markov Chain $X(t)$ on the set $V$
of vertices, whose jump rates are given, when $x$ and $y$ are
neighbours on the graph, by:
$$ w_{x,y} = \nu e^{-\beta ((1-a)E_x -aE_y)} $$
and $ w_{x,y} = 0$ if $x$ and $y$ are not neighbours.

Here $\nu>0 $ is a time-scale parameter, often set to 1, and $a\in [0,1]$ is a symmetry index.

One can also write
$$ w_{x,y} = \nu \tau(x)^{-(1-a)} \tau(y)^a.$$

These jump rates satisfy the detailed balance equation:
$$ \sum_{y\in V} w_{x,y} \tau(x) = \sum_{y\in V} w_{y,x} \tau(y).$$

So that the measure $\tau$ is reversible for the Markov Chain $X$,
whatever the value of the parameter a.

Notice nevertheless that the case where $a=0$ is simpler. Then the
jump rate $w_{x,y}= \frac{\nu}{\tau(x)}$ depends only on the
random landscape at site $x$. The Markov Chain $X$ is simply
obtained by a random time-change from the simple random walk on
the graph. This case is always easier to handle, and will be
called the Random Hopping Times dynamics (RHT dynamics). In this
case, $\tau(x)$ is simply the mean time spent at $x$.

\noindent{\bf Remark} The important feature of the random
variables $\tau (x)$ is the fact that they are heavy-tailed:
$$ P(\tau (x) > a)= \frac{1}{a^{\alpha}} $$
in particular their expectation is infinite, when $\alpha<1$.
The whole model could in fact be directly defined in
terms of the $\tau(x)$ rather than from the energies $(-E_x)$, and then one could assume that the $\tau(x)$ are
i.i.d and in the domain of attraction of an $\alpha$-stable law.

\subsection{Aging and two-point functions} \label{subsection 2.3}

\vskip-5mm \hspace{5mm}

The natural question about Bouchaud's Markov Chain is to study its
long-time behaviour, either in the ``quenched" regime, i.e almost
surely in the randomness of the energy landscape, or in the
``annealed" regime, i.e after averaging in this randomness. One
would expect a dynamical phase transition, between the high
temperature phase where $\alpha > 1$ and the low temperature phase
where $\alpha<1$ . A general idea about the low temperature phase
is that the system spends most of its time in very deep traps:
more precisely, that by time t, the chain has explored a large
part of the space and has found traps of depth depending on t. So
that , at age t, with high probability the system sits waiting in
a deep trap, whose depth is t-dependent, for a time thus depending
on t, before being able to get out and find another deep trap. The
beautiful idea put forward by the physics litterature is ``to
think in the two-times plane", i.e to consider the evolution of
the system between two large times, generally denoted $t_w$ (like
waiting time, $t_w$ is the age of the system) and $t_w +t$ (t is
then the duration of the observation of the system), and to let
both $t_w$ and t tend to infinity. The next step proposed by
physicists is to choose appropriate two-point functions, i.e
functions of the evolution of the system in the time interval
($t_w$, $t_w + t$), in order to measure how much the system
forgets its past in this time interval. The simplest such
two-point function is the quenched probability that the system is
in the same state at times $t_w$ and $t_w + t$:
$$ R^{\omega}(t_w, t_w + t) = P^{\omega}(X(t_w)=X(t_w + t)).$$
Here the probability is quenched, i.e conditioned on the random medium, and the superscript $\omega$ denotes this
randomness of medium, i.e the i.i.d collection of energies. One will also consider the annealed version of this
two-point function:
$$ R(t_w, t_w + t) = \langle R^{\omega}(t_w, t_w + t) \rangle$$
where $\langle\ ,\ \rangle$ denotes the expectation w.r.t the
medium. One also considers often another two-point function, i.e
the quenched probability that the system has not jumped at all in
the time interval ($t_w$, $t_w + t$):
$$ \Pi^{\omega}(t_w, t_w + t) = P^{\omega}(X(t_w)=X(t_w + s), \forall s<t)$$
or its annealed conterpart
$$ \Pi(t_w, t_w + t) = \langle \Pi^{\omega}(t_w, t_w + t)\rangle.$$

Finding an ``agime regime" is then proving that when t is too
small as a function of $t_w$ (typically $t=o(t_w^{\gamma})$ for
some ``aging exponent" $\gamma$), then such a two-point function
is close to one, and when it is large enough (typically
$t>>t_w^{\gamma}$) it is close to zero.  Naturally it is even more
desirable to find the limit of these two-point functions for the
critical regime (typically $t=Ct_w^{\gamma}$).

This program is now understood rigorously for a few important
graphs, which we will now review. We will review the case of
$Z^d$, first for $d=1$ after the work of Fontes-Isopi-Newman
([16],[17])and more recently Cerny (see [4],[14]), then the case
of $d=2$ (see [5]) (for $d>2$ see [14]). Then we treat the case of
the complete graph on $M$ points when $M$ tends to infinity. This
is in fact the original Bouchaud model, which was introduced as an
ansatz for the dynamics of the Random Energy Model. We will then
survey the recent results on these dynamics of the REM
([1],[2],[3]) which will be the first results really pertaining to
the topic of spin glass dynamics.

\subsection{Two aging regimes for Bouchaud's model on {\boldmath $Z$}} \label{subsection 2.4}

\vskip-5mm \hspace{5mm}

Bouchaud's model on $Z$ has been first studied by
Fontes-Isopi-Newman, when $a=0$, i.e for the Random Hopping Times
dynamics. To understand aging in dimension $d=1$, it is important
to introduce a limiting object which will play the role of the
random medium:

\begin{definition} The random speed measure $\rho$.

Let $(x_i,v_i)$ be a Poisson Point Process on $R \times (0,\infty)$,
with intensity measure $\alpha v^{-(1+\alpha)} dxdv$.
We define a random measure $\rho$ on $R$ by
$$\rho = \sum_i v_i \delta_{x_i}.$$
\end{definition}

\begin{definition} The FIN (Fontes-Isopi-Newman) singular diffusion $Z(s)$.

Let $W(t)$ be a standard one-dimensional Brownian Motion, and
l($t$,$y$) its local time at $y$. Define the random time-change:
$$\phi^{\rho}(t)= \int l(t,y)\rho(dy)$$
and its inverse
$$ \psi^{\rho}(t)= inf(s, \phi^{\rho}(s)=t).$$
Then the FIN singular diffusion is $Z(s)= W( \psi^{\rho}(t)).$
\end{definition}

Notice that the Random speed measure and the FIN singular
diffusion are entirely independent of the symmetry parameter $a$,
but depend only on the temperature parameter $\alpha$.

Then the following (annealed) aging result has been proved in [16]
for $a=0$, and in [4] for general a's.

\begin{theorem}
For any $\alpha<1$ and any $a\in [0,1]$, the following limit exists
\[
\lim_{t_w \to \infty} R(t_w,(1+\theta)t_w)= f(\theta).
\]
Moreover the function f can be computed using the singular
diffusion $Z$:
\[
f(\theta)= \langle P(Z(1+\theta)=Z(1)) \rangle.
\]
\end{theorem}

This result shows that the two-point function $R$ exhibits an
aging regime, $t=\theta t_w$, independently of $a$. To be able to
feel the influence of $a$, one should use the other (annealed)
two-point function $\Pi$, which exhibits another aging regime (see
[4]). Let us introduce some notation: Denote by $F$ the annealed
distribution function of the important r.v $\rho(Z(1)$
$$F(u)=\langle P(\rho(Z(1)\leq u) \rangle.$$
Here the brackets denote the average w.r.t the environment i.e the randomness of the measure $\rho$.
Let $g_a$ be the Laplace transform of the r.v $\tau(0)^a$
$$g_a(\lambda) = E[ e^{-\lambda \tau(0)^a}]$$
and $C$ the constant given by $C=2^{a-1}[E(\tau(0)^{-2a})]^{1-a}.$

\begin{theorem}
For any $\alpha<1$ and any $a\in [0,1]$, the following limit exists
\[
\lim_{t_w \to \infty} \Pi(t_w,t_w + \theta t_w^{\gamma})= q_a(\theta)
\]
where $\gamma=\frac{1-a}{1+\alpha}$ and this limit can be computed
explicitly
\[
q_a(\theta) = \int_0^\infty g_a^2(C\theta u^{a-1})\,dF(u).
\]
In particular, when $a=0$,
\[
q_0(\theta) = \int_0^\infty e^{-\frac{\theta}{u}}\,dF(u).
\]
\end{theorem}

These two results show that the Markov Chain essentially succeeds
in leaving the site it has reached at age $t_w$ only after a time
$t$ of the order of $t_w^{\gamma}$, but that for time scales
between $t_w^{\gamma}$ and  $t_w$ it will jump out of the site
reached at age $t_w$ but will not find an other trap deep enough
and so will be attracted back to the trap it reached at age $t_w$.

\subsection{Aging for Bouchaud's model on {\boldmath $Z^2$}, the RHT case} \label{subsection 2.5}

\vskip-5mm \hspace{5mm}

In dimension 2, i.e on $Z^2$, the only case studied is the case where a=0, the RHT case.
An aging regime is exibited (in [5])for both quenched two-point functions

\begin{theorem}
The following limits exists almost surely in the environment $\omega$,
$$ \lim_{t_w \to \infty} R^{\omega}(t_w,t_w+ \theta t_w)= h(\theta),$$
$$ \lim_{t_w \to \infty} \Pi^{\omega}(t_w,t_w+ \theta \frac{t_w}{lnt_w})= k(\theta).$$

The functions h and k satisfy
$$ \lim_{\theta \to 0} h(\theta)= \lim_{\theta \to 0} k(\theta)=1$$
and
$$ \lim_{\theta \to \infty} h(\theta)= \lim_{\theta \to \infty} k(\theta)=0.$$
\end{theorem}
In fact the functions $h$ and $k$ can be computed explicitly
easily using arcsine laws for stable processes. For instance
$$ h(\theta)= \frac{\sin\alpha\pi}{\pi}\int_0^{1+\theta} u^{\alpha -1}(1-u)^{-\alpha}\,du.$$

Here again there is a difference between the subaging regime ( i.e $t= \theta \frac{t_w}{ln(t_w)}$for $\Pi$))and
aging regime(i.e $ t=\theta t_w$ for R), but much slighter than in dimension 1. Indeed it is naturally much more
difficult to visit a trap again after leaving it.

There is another important difference between $d$=1 and $d$=2,
noticed by [16], namely there is localisation in dimension 1 and
not in dimension 2. More precisely, in $d$=1
$$ \limsup_{t \to \infty} \sup_{x \in Z}\langle P^{\omega}(X(t)=x) \rangle >0.$$
And this property is wrong in dimension 2.

\subsection{Aging for Bouchaud's model on a large complete graph} \label{subsection 2.6}

\vskip-5mm \hspace{5mm}

Consider now the case where $G$ is the complete graph on $M$
points. We will also study here only the case the RHT case, where
$a$=0. We consider Bouchaud's Markov Chain $X$ on $G$, started
from the uniform measure. Then it is easy to see that the times of
jump form a renewal process and that the two-point function
$\Pi^{\omega}_M$ is the solution of a renewal equation. This
renewal process converges, when $M$ tends to $\infty$, to a
heavy-tailed renewal process we now introduce.

Let
$$F_{\infty}(t) = 1-\alpha \int_1^{\infty} e^{-\frac{t}{x}}x^{-(1+\alpha)}\,dy.$$
Consider $\Pi_{\infty}(t_w,t_w+t)$ the unique solution of the renewal equation
$$ \Pi_{\infty}(t_w,t_w+t)=1- F_{\infty}(t_w + t) + \int_0^t \Pi_{\infty}(t_w-u,t_w-u+t)\,dF_{\infty}(u).$$

\begin{theorem}{\rm (see [10], [3])}
Almost surely in the environment, for all $t_w$ and t,
$\Pi^{\omega}_M(t_w,t_w+t)$ converges to
$\Pi_{\infty}(t_w,t_w+t).$
\end{theorem}

It is easy to see that the limiting two-point function $\Pi_{\infty}(t_w,t_w+t)$ shows aging.
Let
$$ H(\theta)= \frac{1}{\pi cosec(\frac{\pi}{\alpha})} \int_{\theta}^{\infty} \frac{1}{(1+x)x^{\alpha}\,dx}.$$
Obviously $H(\theta)\sim 1- C_{\alpha}\theta^{1 - \alpha}$ when $\theta$ tends to 0, and
$ H(\theta)\sim C'_{\alpha}\theta^{ - \alpha}$ when $\theta$ tends to $\infty$ .
\begin{theorem}
 $$\lim_{t_w \to \infty}\Pi_{\infty}(t_w,t_w+\theta t_w)= H(\theta ).$$
So that
$$\lim_{t_w \to \infty} \lim_{M \to \infty } \Pi^{\omega}_M (t_w,t_w+\theta t_w)= H(\theta ).$$
\end{theorem}

We will see in the next section how this simple result gives an approximation (in a weak sense) for the much more
complex problem of aging for the RHT dynamics for the REM.

\section{Aging for the Random Energy Model} \label{section 3}

\vskip-5mm \hspace{5mm}

We report here on the joint work with Anton Bovier and Veronique Gayrard, see [1],[2], [3].
This work uses heavily the general analysis of metastability for disordered mean-field models in [11], [12].
Let us first, following Derrida, define the REM, often called the simplest model of a spin glass. A spin configuration $\sigma$ is a vertex of the hypercube $ S_N = \{-1,1\}^N $. As in Bouchaud model on the graph $ S_N$ we consider a collection of i.i.d random variables $(E_{\sigma})_{\sigma \in {S_N} }$ indexed by the vertices of $S_N$.
But we will here assume  that the distribution of the $E_{\sigma}$ 's is standard Gaussian.
We then define the Gibbs measure $\mu_{\beta,N}$ on $S_N$ by setting:
$$ \mu_{\beta,N}(\sigma) = \frac{e^{\beta \sqrt{N} E_{\sigma} }}{Z_{\beta,N} } $$
where $Z_{\beta,N}$ is the normalizing partition function. The
statics of this model are well understood (see[15], and [13]). It
is well known that the REM exhibits a static phase transition at
$\beta_c = \sqrt{2 \ln2}$. For $\beta > \beta_c$ the Gibbs measure
gives ,asymptotically when $N$ tends to $\infty$, positive mass to
the configurations $\sigma$ where the extreme values of the order
statistics of the i.i.d $N(0,1)$ sample $(E_{\sigma} )_{\sigma \in
{S_N} } $ are reached, i.e if we order the spin configurations
according to the magnitude of their (-) energies:
$$ E_{\sigma^{(1)}} \geq  E_{\sigma^{(2)}} \geq ... \geq E_{\sigma^{(2^N)}}.$$

Then for any fixed $k$, the mass $\mu_{\beta,N}(\sigma^{(k)}$ will
converge to some positive random variable. In fact the whole
collection of masses $\mu_{\beta,N}(\sigma^{(k)})$ will converge
to a point process, called Ruelle's point process. Consider for
any $E$, the set  (which we have called the top in [2] and [3],
and should really be called the bottom) of configurations with
energies below a certain threshold $u_N (E)$.

$$T_N(E) = \{ \sigma \in S_N, E_{\sigma} \geq u_N (E) \} .$$
We will choose here the natural threshold for extreme of standard Gaussian i.i.d rv's, i.e
$$ u_N (E)= \sqrt{2N \ln2} + \frac{E}{\sqrt{2N \ln2}} - \frac{\ln(N \ln2) + \ln(4 \pi)}{2 \sqrt{2N \ln2}}.$$
Now we define  discrete-time dynamics on $S_N$ by the transition probabilities
\begin{equation}
p_N (\sigma, \eta)=
\begin{cases}
\frac{1}{N} e^{- \beta \sqrt{N} E_{\sigma}^+ }, &\text{if $\sigma$ and $\eta$ differ by a spin-flip;} \\
1 - e^{- \beta \sqrt{N} E_{\sigma}^+ }, &\text{if $\sigma = \eta$;} \\
0, &\text{otherwise.}
\end{cases}
\end{equation}
Notice we have here truncated the negative values of the
$E_{\sigma}$'s , this truncation is technical and irrelevant. We
could truncate much less drastically. Anyway the states with very
negative $E$'s wont be seen on the time scales we are interested
in.

Then the idea defended by Bouchaud is that the motion of these REM
dynamics when seen only on the deepest traps $T_N(E)$ should be
close to the dynamics of the Bouchaud model on the complete graph
for large $M$. This is true only to some extent. It is true that,
if one conditions on the size of the top $T_N (E)$ to be $M$, then
the sequence of visited points
 in $T_N (E)$ has asymptotically, when $N$ tends to $\infty$, to the standard random walk on the complete graph with $M$ points. Nevertheless Bouchaud's picture would be completely correct if the process observed on the top would really be Markovian, which is not the case, due to a lack of time scales separation between the top and its complement.
Nevertheless it is remarkable that in a weak asymptotic form
Bouchaud's prediction about aging is correct. Let us consider the
following natural two-point function:
$$ \Pi_N (n,m) = \frac{1}{|T_N (E)|} \sum_{\sigma \in T_N (E)} \Pi_{\sigma}(n,m)$$
where $\Pi_{\sigma} (n,m)$ is the quenched probability that the
process starting at time 0 in state $\sigma$ does not jump during
the time interval $(n, n+m)$ from one state in the top $T_N (E)$
to another such state.
\begin{theorem} Let $ \beta > \beta_c $
Then there is a sequence $ c_{N,E} \sim e^{\beta \sqrt{N} u_N (E)} $, such that for any $\epsilon > 0$
$$ \lim_{t_w \to \infty} \lim_{E \to -\infty}
 \lim_{N \to \infty } P (|\frac{\Pi_N (c_{N,E} t_w, c_{N,E} (t_w + t))}{H(\frac{t}{t_w}}) - 1| \geq \epsilon )= 0.$$
\end{theorem}

\begin{remark}\rm
The rescaling of the time by the factor $c_{N,E}$ shows that Bouchaud's trap model is a good approximation of the REM dynamics for the very large time asymptotics, on the last time scale before equilibrium is reached. This is to be contrasted with the other model of glassy dynamics for spin glasses, as advocated by Parisi, Mezard or Cugliandolo-Kurchan where the infinite volume limit is taken before the large time limit. We will see in the next section an example of such a very different aging phenomenon in a much shorter time-scale.
\end{remark}

\begin{remark}\rm
The REM is the first model on which this study of aging on very long times scales (activated dynamics)has been rigorously achieved.
This phenomenon should be present in many more models, like the Generalized Random Energy Model, which is the next achievable goal, or even in much harder problems like the p-spin models for large enough p's. The tools developped in [11][12] should be of prime relevance.
\end{remark}

\section{Aging for the spherical Sherrington-Kirkpatrick model} \label{section 4}

\vskip-5mm \hspace{5mm}

Studying spin glass dynamics for the Sherrington Kirkpatrick model
might seem premature, since statics are notoriously far from fully
understood, as opposed to the REM. Nevertheless, following
Sompolinski and Zippelius, a mathematical study of the Langevin
dynamics has been undertaken in the recent years jointly with
A.Guionnet (see [7],[8],[19]). The output of this line of research
has been to prove convergence and large deviation results for the
empirical measure on path space as well as averaged and quenched
propagation of chaos. The same problem has been solved by $M$.
Grunwald for discrete spins and Glauber dynamics, see [18]. The
law of the limiting dynamics (the self consistent single spin
dynamics) is characterized in various equivalent ways, from a
variational problem to a non-Markovian implicit stochastic
differential equation, none of which being yet amenable, for the
moment, to a serious understanding.

The Sherrington-Kirkpatrick Hamiltonian is given by
$$H^N_J(x)=\frac{1}{\sqrt{N}} \sum_{i,j=1}^N J_{ij} x_i x_j$$
with a $N \times N$ random matrix ${\bf J}=(J_{ij})_{1 \leq i,j \leq N}$
of centered i.i.d standard Gaussian random variables.
The Langevin dynamics for
this model are described by the
stochastic differential system  :

\begin{equation}\label{interaction}
dx^j_t=dB_t^j- U'(x^j_t) dt - \frac{\beta}{\sqrt N} \sum_{1\leq
i\leq N} J_{ji}x^i_t dt,
\end{equation}
where $B$ is a $N$-dimensional Brownian motion, and U a smooth
potential growing fast enough to infinity. It was proved in
[8],[19] that, for any time $T>0$,
 the empirical measure on path space
$$\mu _N:=\frac{1}{N} \sum_{i=1}^N \delta_{x^i_{[0,T]}}$$
converges almost surely towards a non Markovian limit law
$Q_{\mu_0}^T$, when the initial condition is a ``deep quench", i.e
when $x_0=\{x^j_0,1\leq j\leq N\}$ are i.i.d. $Q_{\mu_0}^T$ is
called the self consistent single spin dynamics in the physics
litterature. It is the law of a self-consistent non Markovian
process, which is very hard to study. One expects that the long
time behaviour of this process shows an interesting dynamical
phase transition, and in particular exhibits aging. A consequence
of the convergence stated above is that the limit of the empirical
covariance exists, and is simply the autocovariance of the law
$Q_{\mu_0}^T$ :
$$ C(t_w, t_w + t) := \lim_{N \to \infty} \frac{1}{N} \sum_{i=1}^N x_{t_w}^i x_{t_w +t}^i = \int x_{t_w}x_{t_w +t}\,dQ_{\mu_0}^T(x).$$

Unfortunately it is not possible to find a simple, autonomous equation satisfied by C, or even by C and the so-called response function R. This is a very hard open problem.
But in the physics litterature (mainly in the work of Cugliandolo and Kurchan) one find that the same program is tractable and gives a very rich picture of aging for a large class of models, i.e the spherical p-spin models.

We report here on the joint work with A.Dembo and A.Guionnet [6], on the simplest of such models, i.e the Spherical SK, or spherical p-spin model with p=2. The general case of $p>2$ is harder and will be our next step in the near future.
In this work, we study the Langevin dynamics for a  spherical version
of the Sherrington Kirkpatrick (SSK) spin glass model.

More precisely, we shall consider the
following stochastic differential system
\begin{equation}
du^i_{t}= \beta \sum_{j=1}^N J_{ij}
u^j_t dt - f'(\frac{1}{N}\sum_{j=1}^N (u^j_t)^2)
u^i_t dt  + dW^i_t
\end{equation}
where $f'$ is a uniformly Lipschitz, bounded below
function on $R^+$ such that $f(x)/x \to \infty$ as $x \to \infty$ and
$(W^i)_{1\le i\le N}$ is an $N$-dimensional Brownian motion,
independent of $\{ J_{i,j} \}$ and of the initial data  $\{ u^i_0 \}$.

The term containing $f$ is a Lagrange
multiplier used to implement a ``soft'' spherical constraint.

Here again the empirical covariance admits a limit
$$ C(t_w, t_w + t) := \lim_{N \to \infty} \frac{1}{N} \sum_{i=1}^N x_{t_w}^i x_{t_w +t}^i .$$
But now, as opposed to the true SK model, this limiting two point function is easily computable from an autonomous renewal equation.
Indeed the induced rotational symmetry of the spherical model
reduces the dynamics in question to an $N$-dimensional coupled system
of Ornstein-Uhlenbeck processes whose random drift parameters are
the eigenvalues of a GOE random matrices.
\begin{theorem} There exists a critical $\beta_c$ such that:
When starting from i.i.d initial conditions:
If $\beta <\beta_c$, then $$C(t_w,t_w+t) \leq C_{\beta} \exp(-\delta_\beta |t-s|)$$ for some
$\delta_\beta>0$, $C_\beta <\infty$ and all $(t_w,t)$.

If $\beta=\beta_c$, then $C(t_w,t_w+t) \to 0$ as $t \to \infty$. If
$t_w/t$ is bounded, then the  decay is polynomial $t^{-1/2}$,
and otherwise it behaves like $\frac{t_w^{1/2}}{t_w + t}$.

If $\beta >\beta_c$ then the following limit exists $$ \lim_{t_w
\to \infty} C(t_w, t_w+\theta t_w) = f(\theta).$$

Moreover if $t \gg t_w \gg 1$, then $C(t_w,t_w+t) \frac{t}{t_w}^{3/4}$ is
bounded away from zero and infinity. In particular, the convergence of
$C(t_w,t_w+t)$ to zero occurs if and only if $\frac{t}{t_w} \to \infty$.
\end{theorem}

In these much shorter time scales than for the former sections, the aging phenomenon we exhibit here
is quite different than the one shown in the REM. Here the system has no time to cross any barrier, or to explore and find deep wells,it simply goes down one well, in very high dimension. In some of these very many directions (corresponding to the top eigenvectors of the random matrix {\bf J}), the curvature of the well is so weak that the corresponding coordinates of the system are not tightly bound and are very slow to equilibrate. These very slow  components are responsible for the aging phenomenon here.

\label{lastpage}

\end{document}